\documentclass{article}
\usepackage[utf8]{inputenc}
\usepackage{graphicx}
 
\usepackage{graphicx}
\usepackage{tabularx}
\usepackage{amsmath}
\usepackage{xcolor}
\usepackage{booktabs}
\usepackage{color, soul}
\usepackage{amssymb,latexsym}
\usepackage{amsmath,amsthm}
\usepackage{enumerate}
\usepackage{caption}
\usepackage{subcaption}
\usepackage[labelfont=bf]{caption}
\usepackage[labelsep=period]{caption}
\usepackage[mathscr]{eucal}
\usepackage{float}
\floatstyle{plaintop}
\restylefloat{table}
\usepackage[nobysame]{amsrefs}

\restylefloat{table}
\usepackage{placeins}
\usepackage{graphics} 
\usepackage{graphicx}
\graphicspath{ {./images/} }
\usepackage{float}
\title{Natural vibrations of circular nano-arches of piecewise constant thickness}
\author{Jaan Lellep and Shahid Mubasshar}
\date{}

\begin{document}
\maketitle

\section*{Abstract}
The free vibrations of elastic circular arches made of a nano-material are considered. A method of determination of eigenfrequencies of nano-arches weakened with stable cracks is developed making use of the concept of the massless spring and Eringen's nonlocal theory of elasticity. The aim of the paper is to evaluate the sensitivity of eigenfrequencies on the geometrical and physical parameters of the nano-arch.

\section{Introduction}
During the last decades, a lot of attention is paid to the investigation of the structural elements on the nano-size level. It is evident that the classical theory of elasticity is not able to account for the size effects of nano-particles. However, a suitable tool for this purpose presents the nonlocal theory of elasticity developed by Eringen \cite{eringen2002nonlocal}, Eringen and Edelen \cite{eringen1972nonlocal}. Originally Eringen and Edelen \cite{eringen1972nonlocal} applied the nonlocal constitutive model for the investigation of surface waves and screw dislocations in solids. Later the nonlocal beam theory was used for the investigation of the bending problems by Reddy \cite{reddy2007nonlocal}, Thai \cite{thai2017review}, Li and Wang \cite{li2008introduction}, Li et al. \cite{li2015comments}, Reddy and Pang \cite{reddy2008nonlocal}, Civalek und Demir \cite{civalek2011bending} by making use of the Euler-Bernoulli beam model, also by Roque et al. \cite{roque2011analysis}, Wang, Zhang and He \cite{wang2007vibration} in the case of the Timoshenko beam model. Buckling and vibration of nano-beams were studied by many investigators, including Ansari and Sahmani \cite{ansari2011bending}, Roostai and Haghpanahi \cite{roostai2014vibration}, Sahmani and Ansari \cite{sahmani2011nonlocal}, Behera and Chakraverty \cite{behera2015application}, Aydogdu \cite{aydogdu2009general}, Thai \cites{thai2017review, thai2012nonlocal}, Thai, Vo, Nguyen and Kim \cite{thai2017review}, Murmu and Pradhan \cite{murmu2009small}, Wang, Zhang and He \cite{wang2007vibration}, Ganapathi and Polit \cite{ganapathi2018nonlocal}, Ba\u{g}datli \cite{baugdatli2015non}, Wang, Zhang, Ramesh and Kitipornchai \cite{wang2006buckling}, Wang, Zhang, Challamel and Duan \cite{wang2017boundary}, Wang and Arash \cite{wang2014review} and others. Buckling of multi-step non-uniform beams is investigated by Li \cite{li2001buckling}.  The Ritz method was adopted for the solution of bending, buckling and vibration in the case of nano-beams by Ghannadpour, Mohammadi and Fazilati \cite{ghannadpour2013bending}. It appeared that together with the Ritz method another effective method of solution of problems of this kind is the differential quadrature method (see Shu \cite{shu2000differential}, Pradhan and Kumar \cite{pradhan2011vibration}). The same numerical procedure was employed by Hossain and Lellep \cites{hossain2021mode, hossain2022analysis} as well when studying the natural frequencies of nano-beams with cracks. Note that in these studies, the effect of the rotatory inertia of an element was taken into account in contrast to the classical approach. In the published literature, one can find a limited number of papers concerning the vibration and stability of nano-beams, nano-plates and nano-rings with cracks. An analytical method is developed by Loghmani and Yazdi \cite{loghmani2018analytical} for the determination of natural frequencies of cracked nano-beams with stepped cross-sections. Nano-beams with different end conditions are treated. The case of a cantilever nano-beam with a buckyball at the free end is studied separately in \cite{hosseini2018review}. Curved beams and segments of rings with cracks are the subjects of investigations by Karaagac, Öztürk and Sabuncu \cites{karaagac2011crack, karaagac2009free}, also by Cerri and Ruta \cite{cerri2004detection} and Cerri, Dilena and Ruta \cite{cerri2008vibration}, also Krawczuk et al. \cite{krawczuk1997dynamics}. In-plane free vibrations of circular arches with defects are investigated by Viola, Artioli, and Dilena \cite{viola2005analytical} and Viola and Tornabene \cite{viola2009free}, Chondros et al. \cite{chondros1998continuous}, also by Mazanoglu and Sabuncu \cite{mazanoglu2009vibration} and Mazanoglu et al. \cite{mazanoglu2010flexural}. The influence of defects on the natural frequency of nano-rings was assessed by Wang and Duan \cite{wang2008free} and Moosavi et al. \cite{moosavi2011vibration}. Here the defects are modelled with the help of hinges having rotational restraints on vibrations of nano-rings.
In the present paper, the natural vibration is investigated under the assumption that the nano-arches have step-wise varying cross-sections and that the arches are weakened with crack-like defects. The cracks are assumed to be stationary cracks. No attention will be paid to the extension of these defects. However, the additional compliance induced by cracks is taken into account with the help of the local compliance matrix.

\section{Problem formulation}
\begin{figure}[H]
\centerline{\includegraphics[width=14cm,height=10cm]{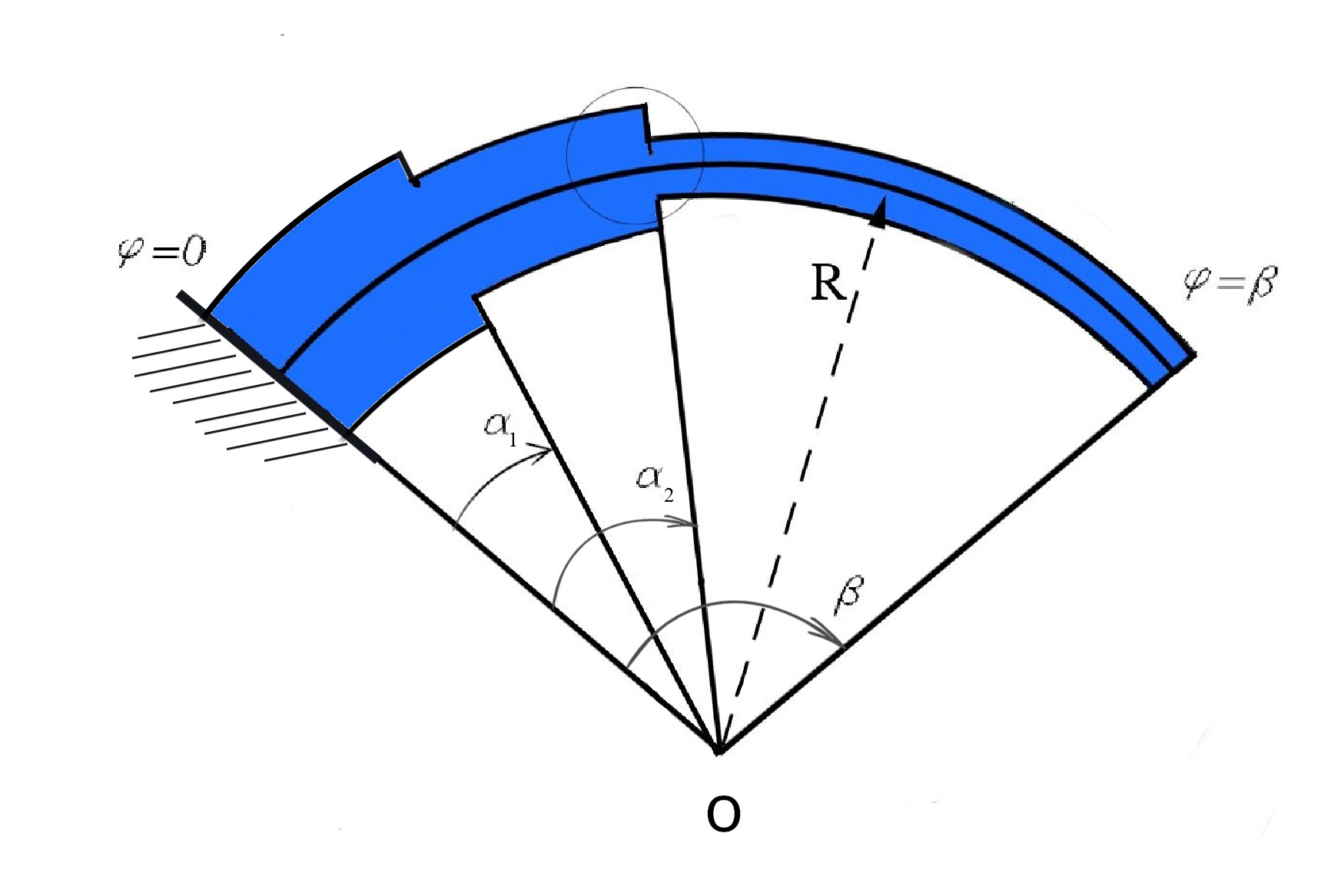}}
\caption{The geometry of the cantilever nano-arch with defects.}
\label{Fig.1(a)}
\centering
\end{figure}

Let us consider the dynamic behaviour of a nano-arch or curved nano-beam of piece-wise constant thickness (Fig. \ref{Fig.1(a)})

\begin{figure}[H]
\centerline{\includegraphics[width=8cm,height=6cm]{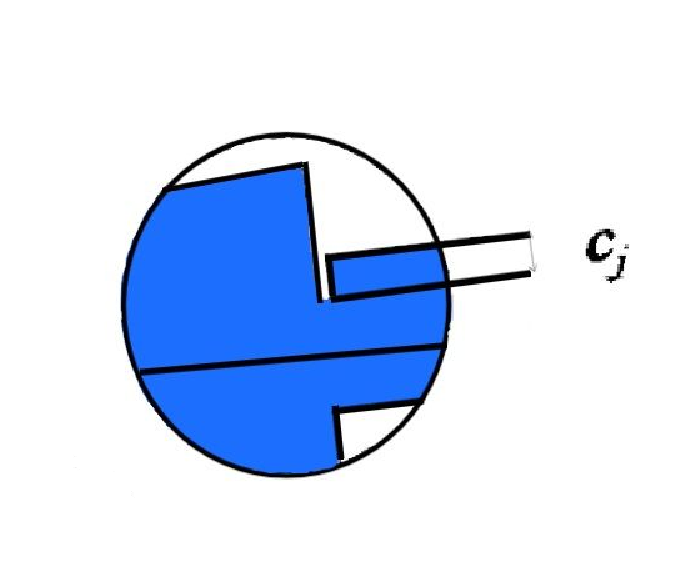}}
\caption{A flaw at the re-entrant corner of the step.}
\label{Fig.1(b)}
\centering
\end{figure}

\begin{equation}\nonumber
\label{eqn:1}
h=
\begin{cases}
    \begin{alignedat}{4}
    &h_{0}\quad&,\quad & \varphi \in [\alpha_{0},\alpha_{1})\\                        
    &h_{1}\quad&,\quad &\varphi \in (\alpha_{1}, \alpha_{2})\\
    & --&--&------\\
    &h_{n}\quad&,\quad&\varphi \in (\alpha_{n}, \alpha_{n+1}],\\
	\end{alignedat}
  \end{cases}
\end{equation}  

where $\alpha_{0}=0$, $\alpha_{n+1}=\beta$. Let $R$ be the radius of the nano-arch. Here $h_{i}$ ($i=0, ..., n$) and $\alpha_{j}$ ($j=0, ..., n+1$) stand for given constants whereas $\varphi$ is the current angle $(0 \leq \varphi \leq \beta)$. In the present study, it is assumed that the nano-arch is clamped at the left-hand end and it is absolutely free at $\varphi=\beta$. It is also assumed that at cross-sections $\varphi=\alpha_{i}$ stable cracks of length $c_{i}$ are situated. The attention will be paid to nano-arches of rectangular cross-sections with thickness $h$ and width $b=const$ only. The aim of the study is to elucidate the sensitivity of eigenfrequencies on the geometrical and physical parameters of the nano-arch and also on the location and length of the defects. For these purposes, the governing equations of the nonlocal theory of elasticity developed by Eringen \cite{eringen2002nonlocal} are used.

\section{Nonlocal material behaviour}

It is widely recognized that the nano-materials subjected to the external loadings behave according to the nonlocal constitutive equations of the theory of elasticity (see Reddy \cite{reddy2007nonlocal}, Eringen \cite{eringen2002nonlocal}, Aydogdu \cite{aydogdu2009general}, Ansari and Sahmani \cite{ansari2011bending}, Loghmani and Yazdi \cite{loghmani2018analytical}, Reddy and Pang \cite{reddy2008nonlocal}, Lellep and Mubasshar \cite{mubasshar2022natural}). It is known that in the classical theory of elasticity, the stress tensor is proportional to the strain tensor at each point of the current body. However, in the nonlocal theory of elasticity the stress at the current point depends on the strain at each point of the body. Probably the simplest nonlocal constitutive equation can be presented as (see Reddy \cite{reddy2007nonlocal}, Thai \cite{thai2012nonlocal}, Lellep and Lenbaum \cites{lellep2022free, lellep2019natural, lellep2018free}
\begin{equation}\nonumber
\label{eqn:7}
\sigma_{ij}-\eta\nabla^{2}\sigma_{ij}=\sigma^c_{ij}.
\end{equation}
Here $\sigma_{ij}$ denotes the stress tensor for a nonlocal theory, $\sigma^c_{ij}$ is the stress tensor in the classical elasticity and $\eta=(e_{0}a)^2$.  Here $a$ is the dimension of the lattice of the material and $e_{0}$ stands for a material constant. Whereas, $\nabla^2$ is the nabla operator. In the present study, it is assumed that $e_{0}$ and $a$ and, thus, also $\eta$ are given constants. Evidently, the determination of $\eta$ may be complicated but this is another task. 

Making use of generalized stresses $M$ (radial bending moment) and $N$ (membrane force in the tangential direction), the nonlocal constitutive equation can be converted into
\begin{equation}
\begin{aligned}
\label{eqn:1a}
M-\eta\nabla^{2}M&=M^c,\\
N-\eta\nabla^{2}N&=N^c,\\
Q-\eta\nabla^{2}Q&=Q^c,
\end{aligned}
\end{equation}

where $M^c$, $N^c$ and $Q^c$ stand for the corresponding
quantities in the classical theory of elasticity (see Reddy \cite{reddy2007nonlocal}, Thai \cites{thai2012nonlocal, thai2017review}, Ansari, Rouhi and Sahmani \cite{ansari2011calibration}, Wang et al. \cite{wang2017boundary}).

It is worthwhile to mention that the total number of equations in the system (\ref{eqn:1a}) coincides with the number of generalized stresses necessary for the correct formulation of the current problem.
In the bending problems, the dominant generalized stresses are the bending moments. In the present case according to (\ref{eqn:1}) one has ($s$ is the length of the element)
\begin{equation}
\label{eqn:2}
M-\eta \frac {\partial^2{M}} {\partial s^2}=-EI_{j}(\frac {\partial^2{W}} {\partial s^2}+\frac{W}{R^2}) 
\end{equation}
for $\varphi \in (\alpha_{j}, \alpha_{j+1})$ ; $j=0, ..., n$. In (\ref{eqn:2}) the substitution 
\begin{equation}\nonumber
\label{eqn:3}
M^c=-EI_{j}(\frac {\partial^2{W}} {\partial s^2}+\frac{W}{R^2}) 
\end{equation}
has been introduced. Here $E$ stands for the Young modulus and $I_{j}$ is the moment of inertia of the arch in the section $(\alpha_{j}, \alpha_{j+1})$ and $ds=Rd\varphi$.
Evidently, in the case of circular arches, the generalized stresses are the bending moment $M$ and the membrane force $N$ (see Soedel \cite{soedel2004}, Lellep and Lenbaum \cites{lellep2022free, lellep2018free}), also the shear force $Q$. Since the latter can be eliminated from the governing equations the shear force will not get any attention in the subsequent analysis. The strain components corresponding $N$ and $M$ are (see Soedel \cite{soedel2004}), the relative extension
\begin{equation}
\label{eqn:3a}
\varepsilon=\frac {1}{R}(U^{'}+W)
\end{equation}
and the curvature
\begin{equation}
\label{eqn:4}
\varkappa=\frac {1}{R^2}(U^{'}+W^{''}).
\end{equation}
In (\ref{eqn:3}) and (\ref{eqn:4}) prims denote the differentiation with respect to the current angle, $U$ and $W$ denote the tangential and transverse displacements of the middle surface, respectively, and $R$ is the radius of the arch. In the classical theory of elasticity (see Soedel \cite{soedel2004}) it is assumed that $\varepsilon=0$ and $U^{'}=-W$. Thus,
\begin{equation}
\label{eqn:5}
\varkappa=-\frac {1} {R}(W+W^{''})
\end{equation}
and 
\begin{equation}
\label{eqn:6}
M^{c}=\frac {-Eh^{3}b} {12R^{2}}(W+W^{''}). 
\end{equation}
The nonlocal constitutive law (\ref{eqn:2}) reads as 
\begin{equation}
\label{eqn:7a}
M-\eta M^{''}=M^c
\end{equation}
and
\begin{equation}
\label{eqn:8}
N-\eta N^{''}=N^c.
\end{equation}

\section{Equilibrium equations}

The equilibrium conditions of an infinitesimal element of the arch can be presented as (see Soedel \cite{soedel2004} and Lellep and Liyvapuu \cite{lellep2016natural})

\begin{equation}
\label{eqn:9}
\begin{aligned}
M^{'}-RQ&=0,\\
  N^{'}+Q+p_{s}R&=\rho h_{j}R\ddot{U},\\
 Q^{'}-N+p_{n}R&=\rho h_{j}R\ddot{W},
 \end{aligned}
\end{equation}

where $j=0, ..., n$ and $\varphi \in (\alpha_{j}, \alpha_{j+1})$. Here $p_{s}$, $U$ and $p_{n}$, $W$ stand for the external loads and displacements in the circumferential and normal direction, respectively. Let $\rho$ be the density of the material and dots denote the differentiation with respect to time $t$ and prims with respect to the length of the arch. Thus, 

\begin{equation}\nonumber
\ddot{U}= \frac {\partial^2{U}} {\partial t^2},\quad \ddot{W}=\frac {\partial^2{W}} {\partial t^2},\quad M^{'}=\frac {\partial{M}} {\partial{ s}}, \quad N^{'}=\frac{\partial{N}} {\partial{s}}.
\end{equation}

Evidently, $ds=Rd\varphi$ in (\ref{eqn:3})-(\ref{eqn:9}). Differentiating the first equation in (\ref{eqn:9}) with respect to $\varphi$ and substituting it into the last one leads to the equation

\begin{equation}
\label{eqn:10}
M^{''}-RN+R^{2}(p-\rho h_{j}\ddot{W})=0,
\end{equation}
which must be satisfied for $\varphi \in (\alpha_{j}, \alpha_{j+1})$. In the case of a curved cantilever nano-beam at the free edge the bending moment and the shear force vanish. Thus,

\begin{equation}
\label{eqn:11}
  M(\beta, t)=0
\end{equation}
and 
\begin{equation}
\label{eqn:12}
  Q(\beta, t)=0.
\end{equation}
However, at the root section at $\varphi=0$ one has
\begin{equation}
\label{eqn:13}
  W(0, t)=0
\end{equation}
and 
\begin{equation}
\label{eqn:14}
  W^{'}(0, t)=0.
\end{equation}

We are considering the natural vibrations of nano-arches. Thus, it is reasonable to assume that $p_{s}=p_{n}=0$ in (\ref{eqn:9}) and $p=0$ in (\ref{eqn:10}). Making use of (\ref{eqn:3}), (\ref{eqn:6}), (\ref{eqn:7a}) and (\ref{eqn:10}), one can define
\begin{equation}
\label{eqn:15}
M=\frac{-1}{R^2(1+\eta)}(EI_{j}(W^{''}+W)-h_{j}\eta \rho R^{2}\ddot{W}) 
\end{equation}

for $\varphi \in (\alpha_{j}, \alpha_{j+1})$; $j=0, ..., n$. In the following, we are looking for the solution of governing equations in the particular case when $M=-RN$. This takes place, for instance, in the case when a cantilevered arch is subjected to the concentrated loading directed towards the center of curvature of the middle line of the arch. It is assumed herein also that the membrane force vanishes. Substituting now (\ref{eqn:15}) in (\ref{eqn:10}) leads to the fourth-order equation
 \begin{equation}
\label{eqn:16}
W^{IV}+2W^{''}+W+\frac{\rho h_{j}R^{4}}{EI_{j}}(\eta \ddot W^{''}-\ddot W)=0,
\end{equation}
which is to be solved for $\varphi \in (\alpha_{j}, \alpha_{j+1})$; $j=0, ..., n$.
\section{Additional compliance due to cracks}

When integrating the governing equations and also (\ref{eqn:16}) one has to account for the continuity of variables $W(\varphi,t)$, $M(\varphi,t)$, $Q(\varphi, t)$ which is the consequence of physical requirements. Moreover, the slope $W^{'}(\varphi, t)$ is also continuous except at the cross-sections with defects. Let us denote
\begin{equation}
\label{eqn:17}
\theta_{j}=W^{'}(\alpha_{j}+0, t)-W^{'}(\alpha_{j}-0, t),
\end{equation}
where $j=1, ..., n$.
The jump of the slope of the deflection of the nano-arch
is coupled with the generalized stresses at the cross-section with defects. Following Dimarogonas \cite{dimarogonas1996vibration}, Anifantis and Dimarogonas \cite{anifantis1983stability} the quantities $\theta_{j}$ will be treated as generalized displacements. In the linear elastic fracture mechanics, the generalized stresses $P_{i}$ and generalized displacements $u_{i}$ $(i=1, ..., 6)$ are coupled as
\begin{equation}
\label{eqn:18}
u_{i}=\frac{\partial U_{s}}{\partial P_{i}},
\end{equation}
where $U_{s}$ stands for the strain energy density. The compliance matrix is defined as
\begin{equation}
\label{eqn:19}
c_{ij}=\frac{\partial u_{i}}{\partial P_{j}}.
\end{equation}
Combining (\ref{eqn:18}) and (\ref{eqn:19}) leads to the equations 
\begin{equation}
\label{eqn:20}
c_{ij}=\frac{\partial^{2} U_{s}}{\partial P_{i} \partial P_{j}},
\end{equation}
where $i=1, ..., 6$; $j=1, ..., n$.
The total energy release rate $G$ is defined as (see Broek \cite{broek2012practical}, Anderson \cite{anderson2004fracture}) 
\begin{equation}
\label{eqn:21}
G=\frac{1-\nu^2}{E}(K_{I}^2+K_{II}^2+\frac{K_{III}^2}{1-\nu}),
\end{equation}
where $E$ is the Young modulus, $\nu$ is the Poisson ratio and 
\begin{equation}
\label{eqn:22}
c_{ij}=\frac{b\;\;\partial^{2}}{\partial P_{i} \partial P_{j}}  \iint_S G \,dS.
\end{equation}
In (\ref{eqn:22}), $c_{ij}$ stand for the elements of the compliance matrix C. The quantities $K_{I}, K_{II}, K_{III}$ stand for the stress intensity factors corresponding to modes I, II and III, respectively (see Broek \cite{broek2012practical}, Anifantis and Dimarogonas \cite{anifantis1983stability}, Alsabbagh, Abuzeid and Dado \cite{alsabbagh2009simplified}). On the other hand, the strain energy release rate can be introduced so that
\begin{equation}
\label{eqn:23}
u_{i}=\frac{b\;\;\partial}{\partial P_{i}} \int_{0}^{c} J \,dc
\end{equation}
and 
\begin{equation}
\label{eqn:24}
c_{ij}=\frac{b\;\;\partial^{2}}{\partial P_{i} \partial P_{j}}  \int_{0}^{c} J \,dc,
\end{equation}

where $J$ is so-called $J-integral$ of Rice (see Rice \cite{rice1968path}). It is worth to emphasize that $dc$ stands in (\ref{eqn:23})-(\ref{eqn:24}) for the extension of the crack and $c_{ij}$ are corresponding compliances due to the crack. Note that the energy release rate $G$ was originally introduced for the characterization of the energy behaviour in the neighbourhood of the crack tip in linear elastic materials.
However, James Rice \cite{rice1968path} was able to extend the energy release rate to non-linear materials. He showed that the non-linear energy release rate can be expressed as a line integral called $J-integral$. The $J-integral$ is evaluated along an arbitrary closed contour $\Gamma$ around the crack tip. It is defined as (see Broek \cite{broek2012practical}, Anderson \cite{anderson2004fracture})
\begin{equation}
\label{eqn:25}
J=\int_{\Gamma} (U_{s}\,dy-T_{i} \frac{\partial u_{i}}{\partial x}\,ds).
\end{equation}
In (\ref{eqn:25}), $T_{i}$ stand for external forces (tractions) and $\Gamma$ is a closed contour followed counterclockwise.
Consider now a particular case when at $\varphi=\alpha_{j}$ a mode I crack is located. If $\theta_{j}$ is treated as a generalized coordinate then the generalized forces associated with $\theta_{j}$ are $M_{j}=M(\alpha_{j}, t)$ and $N_{j}=N(\alpha_{j}, t)$. Thus
\begin{equation}
\label{eqn:26}
\theta_{j}= C_{0j} \overrightarrow{Q} (\alpha_{j}, t),
\end{equation}

where $\overrightarrow{Q}$ is the vector and $M(\alpha_{j}, t), N(\alpha_{j, t})$ and $c_{0j}$ stand for the elements of  corresponding compliance matrix. It is reasonable to introduce the notation

\begin{equation}
\label{eqn:27}
   C_{ij} = \begin{bmatrix}
    c_{11}(\alpha_{j, t}) & c_{12}(\alpha_{j, t})  \\
    c_{21}(\alpha_{j, t}) & c_{22}(\alpha_{j, t})  \\
    \end{bmatrix}.
\end{equation}

Due to symmetry $c_{12}=c_{21}$ at each $\alpha_{j}$. It is worthwhile to mention that the local stiffness matrix $K_{j}$ is reciprocal to the local compliance $c_{j}$ in the one-dimensional case (see Lellep and Sakkov \cite{lellep2006buckling}). An alternative case is studied by Anifantis and Dimarogonas \cite{anifantis1983stability}.
It is known from the linear elastic fracture mechanics that (see Broek \cite{broek2012practical} Anderson \cite{anderson2004fracture}, Broberg \cite{broberg1999cracks})
\begin{equation}
\label{eqn:28}
G=\frac{M^{2}}{2b}\frac{dC}{dc}
\end{equation}

if a beam element is subjected to the bending moment $M$. Here $C$ is the compliance and $c$ stands for the crack length.  In this case the stress intensity factor is
\begin{equation}
\label{eqn:29}
K=\frac{\sigma M}{bh^2} \sqrt{\pi c} F_{1}(\frac{c}{h}).
\end{equation}

In (\ref{eqn:29}), $F_{1}$ stands for so-called shape function which must be determined experimentally, $h$ being the thickness. However, there exist rich databases (see Tada et al. \cite{tada1973stress}) which can be used for the interpolation of shape functions $F_{1}(s)$ and $F_{2}(s)$ (here $s=c/h$). If the beam element is loaded with the moment $M_{j}$ and the tensile force $N_{j}$ then the stress intensity coefficient is calculated as (see Tada et al. \cite{tada1973stress}, Zhou and Huang \cite{zhou2006crack})
\begin{equation}
\label{eqn:30}
K_{j}=\frac{\sqrt{\pi s_{j}}}{b\sqrt{h}}(N_{j}F_{1}(s_{j})+\frac{6 M_{j}}{h}F_{2}(s_{j})),
\end{equation}
where
\begin{equation}
\label{eqn:31}
F_{1}=1.12-0.23s+10.55s^{2}-21.72s^{3}+30.39s^{4}
\end{equation}
and
\begin{equation}
\label{eqn:32}
F_{2}=1.12-1.4s+7.33s^{2}-13.08s^{3}+14.08s^{4}.
\end{equation}
In (\ref{eqn:30}), the quantity $h$ can be interpreted as $h=min(h_{j}, h_{j+1})$.
Another type of the shape function can be presented as (see Tada et al. \cite{tada1973stress}) 
\begin{equation}
\label{eqn:33}
F(s)=\sqrt{\frac{2}{\pi s}{\tan\frac{\pi s}{2}}} \cdot \frac{0.923+0.199(1-\sin\frac{\pi s}{2})^4}{\cos\frac{\pi s}{2}}.
\end{equation}
Making use of the shape functions (\ref{eqn:31})-(\ref{eqn:33}) one can calculate the elements of the compliance matrix (\ref{eqn:27})

\begin{equation}
\label{eqn:34}
\begin{aligned}
c_{11}(\alpha_{j})&=\frac{2}{Eb}  \int_{0}^{c} sF_{1}^2 \,ds,\\
c_{12}(\alpha_{j})&=\frac{2}{Ebh}  \int_{0}^{c} sF_{1}F_{2} \,ds,\\
c_{22}(\alpha_{j})&=\frac{2}{Eb}  \int_{0}^{c} sF_{2}^2 \,ds.
\end{aligned}
\end{equation}
Finally, making use of (\ref{eqn:34}) one can calculate
\begin{equation}
\label{eqn:35}
\theta_{j}=c_{11}(\alpha_{j})M(\alpha_{j}, t)+c_{12}(\alpha_{j})N(\alpha_{j}, t)
\end{equation}
and the corresponding displacement 
\begin{equation}
\label{eqn:36}
\delta_{j}=c_{21}(\alpha_{j})M(\alpha_{j}, t)+c_{22}(\alpha_{j})N(\alpha_{j}, t).
\end{equation}
\section{Solution of the governing equations}
In order to solve the equation (\ref{eqn:16}) in regions $(\alpha_{j}, \alpha_{j+1})$; $j=0, ..., n$ the method of separation of variables is applicable. Assume, thus, that
 \begin{equation}
\label{eqn:37}
  W(\varphi, t)=X_{j}(\varphi) T(t)
\end{equation}
for $\varphi \in(\alpha_{j}, \alpha_{j+1})$; $j=0, ..., n$. The transformation (\ref{eqn:37}) admits to present (\ref{eqn:16}) as the system of equations consisting of
 \begin{equation}
\label{eqn:38}
  \ddot{T}+\omega ^{2}T=0
\end{equation}
and 
\begin{equation}
\label{eqn:39}
X_{j}^{IV}+X_{j}^{''}(2+A_{j}\eta)+X_{j}(1-A_{j})=0,
\end{equation}
where
\begin{equation}
\label{eqn:40}
A_{j}=\frac{\omega ^{2} \rho h_{j} R^{4}}{EI_{j}}; j=0, ..., n.
\end{equation}
The characteristic equation for the linear fourth-order equation (\ref{eqn:39}) is
\begin{equation}
\label{eqn:41}
\lambda_{j}^{4}+(2+A_{j}\eta)\lambda_{j}^{2}+1-A_{j}=0; j=0, ..., n.
\end{equation}
It immediately follows from (\ref{eqn:41}) that 
\begin{equation}
\label{eqn:42}
\lambda_{j}^{2}=\frac{1}{2}(-2-A_{j}\eta)\pm\frac{1}{2}\sqrt{(2+A_{j}\eta)^{2}-4(1-A_{j})}
\end{equation}
and the general solution of (\ref{eqn:39}) has the form 
\begin{equation}
\label{eqn:43}
X_{j}=C_{1j} \cosh\mu_{j} \varphi+C_{2j} \sinh\mu_{j} \varphi+C_{3j} \cos\nu_{j} \varphi+C_{4j} \sin\nu_{j}\varphi,
\end{equation}
where $j=0, ..., n$ and
\begin{equation}
\label{eqn:44}
 \mu_{j}=\sqrt{-1-\frac{\eta}{2}A_{j}+B_{j}} 
\end{equation}
and
\begin{equation}
\label{eqn:45}
 \nu_{j}=\sqrt{1+\frac{\eta}{2}A_{j}+B_{j}}, 
\end{equation}
where the notation 
\begin{equation}
\label{eqn:46}
 B_{j}=\frac{1}{2}\sqrt{A_{j}^2\eta^2+4A_{j}(1+\eta)}
\end{equation}
is introduced. Note that the solution of the equation (\ref{eqn:38}) depending on time can be presented as 
\begin{equation}
\label{eqn:47}
T=\sin (\omega t),
\end{equation}
$\omega$ being the natural frequency of the nano-arch. Evidently, now
\begin{equation}
\label{eqn:48}
W(\varphi, 0)=0,\quad  \dot W(\varphi, 0)=\omega X_{j}(\varphi)
\end{equation}
for $\varphi \in(\alpha_{j}, \alpha_{j+1})$; $j=0, ..., n$. It is presumed herewith that the initial conditions are given in the form (\ref{eqn:48}). Consider now the boundary conditions. At the root section $\varphi=0$ according to (\ref{eqn:13}), (\ref{eqn:14}) and (\ref{eqn:37}), one has
\begin{equation}
\label{eqn:49}
X_{0}(0)=0,\quad  X_{0}^{'}(0)=0.
\end{equation}
The boundary conditions corresponding to the free edge are presented by (\ref{eqn:11}) and (\ref{eqn:12}). Making use of (\ref{eqn:15}) and (\ref{eqn:37}) one can check that the requirements at the free edge are satisfied if
\begin{equation}
\label{eqn:50}
X_{n}^{''}(\beta)+X_{n}(\beta)(1+A_{n})=0
\end{equation}
and 
\begin{equation}
\label{eqn:51}
X_{n}^{'''}(\beta)+X_{n}^{'}(\beta)(1+A_{n})=0.
\end{equation}
\section{Intermediate jump conditions}
It was stated above that due to cracks the nano-arch has additional flexibility. The additional compliance can be calculated by making use of the relations (\ref{eqn:26})-(\ref{eqn:36}). According to this concept, the slope $W^{'}$ has finite jumps $[W^{'}(\alpha_{j})]$ at $\varphi=\alpha_{j}$ $(j=1, ..., n$). Thus
\begin{equation}
\label{eqn:52}
W^{'}(\alpha_{j}+0, t)=W^{'}(\alpha_{j}-0, t)+\theta_{j},
\end{equation}
where $\theta_{j}$ is evaluated by (\ref{eqn:26}), (\ref{eqn:27}). Making use of (\ref{eqn:15}), (\ref{eqn:37}), one can present 
\begin{equation}
\label{eqn:53}
M(\alpha_{j}+0, t)=\frac{-EI_{j}}{(1+\eta)R^2}(X_{j}^{''}(\alpha_{j})+(1+A_{j}\eta)X_{j}(\alpha_{j})T(t).
\end{equation}

In (\ref{eqn:52}), (\ref{eqn:53}) one has to distinguish the left-hand and right-hand limits, respectively. Here the notation
\begin{equation}
\label{eqn:54}
g(\alpha\pm0)=\lim_{x\to\alpha\pm0} g(x)
\end{equation}
is used. It follows from (\ref{eqn:17}), (\ref{eqn:26})-(\ref{eqn:35}) that 
\begin{equation}
\label{eqn:55}
\theta_{j}=(c_{11}(\alpha_{j})-\frac{c_{12}(\alpha_{j})}{R})M(\alpha_{j},t),
\end{equation}
where the compliances $c_{11}(\alpha_{j})$ and $c_{12}(\alpha_{j})$ are evaluated by (\ref{eqn:34}). It is evident from the physical considerations that the stress components $M(\varphi, t)$ and $Q(\varphi, t)$ together with the displacement $W(\varphi, t)$ must be continuous at each time instant at each $\varphi\in (0, \beta)$. Thus, the set of continuity and jump conditions is 
\begin{equation}
[W(\alpha_{j}, t)]=0, \quad [W^{'}(\alpha_{j}, t)]=\theta_{j}, \quad [M(\alpha_{j}, t)]=0, \quad [Q(\alpha_{j}, t)]=0,
\end{equation}
$j=1, ..., n.$
\section{Numerical results and conclusions}
\begin{figure}[H]
\centerline{\includegraphics[width=14cm,height=8cm]{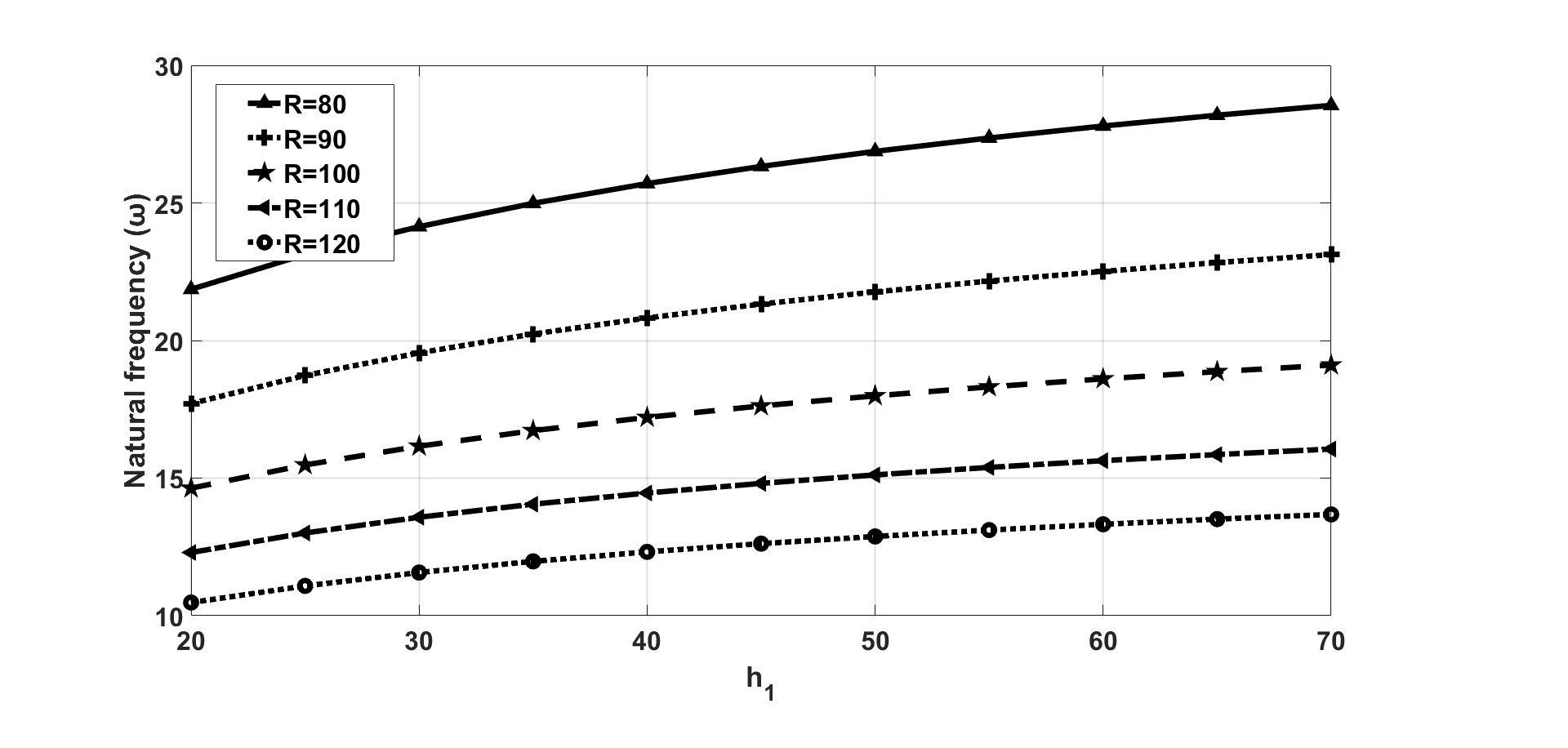}}
\caption{Natural frequency versus thickness of the nano-arch for $\alpha=0.8$.}
\label{Fig.3}
\centering
\end{figure}

\begin{figure}[H]
\centerline{\includegraphics[width=14cm,height=8cm]{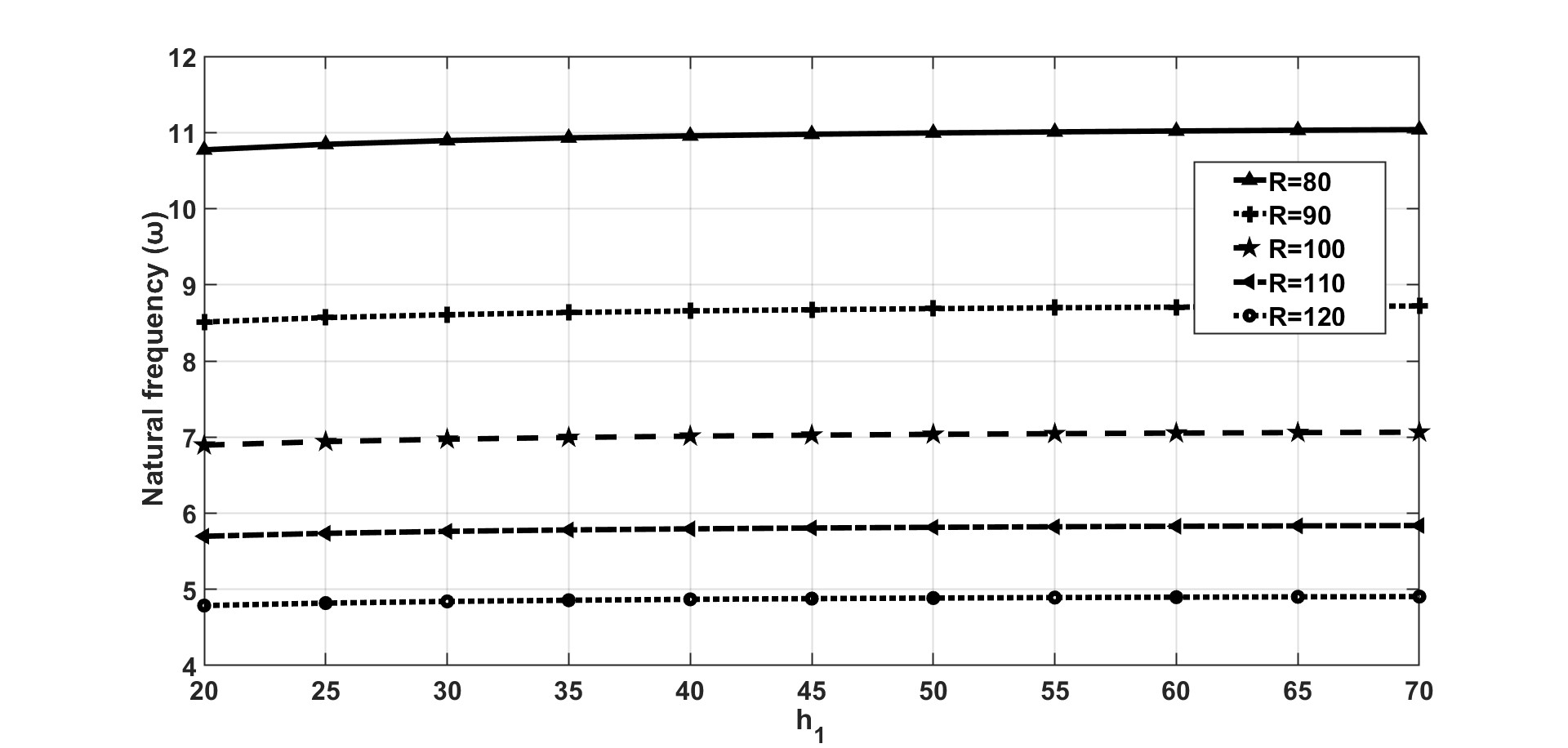}}
\caption{Natural frequency versus thickness of the nano-arch for $\alpha=0.6$.}
\label{Fig.4}
\centering
\end{figure}

\begin{figure}[H]
\centerline{\includegraphics[width=14cm,height=8cm]{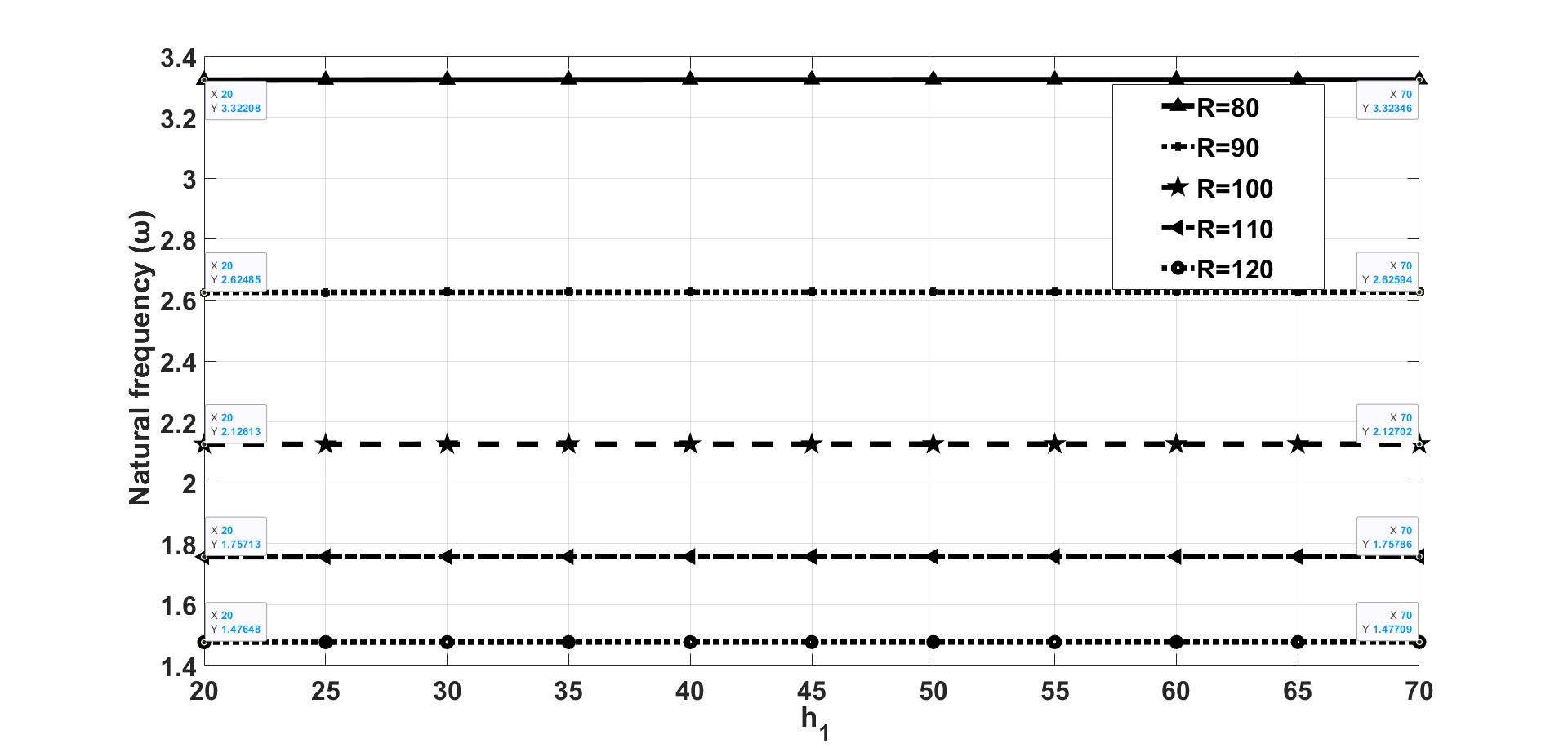}}
\caption{Natural frequency versus thickness of the nano-arch for $\alpha=0.4$.}
\label{Fig.5}
\centering
\end{figure}

\begin{figure}[H]
\centerline{\includegraphics[width=14cm,height=8cm]{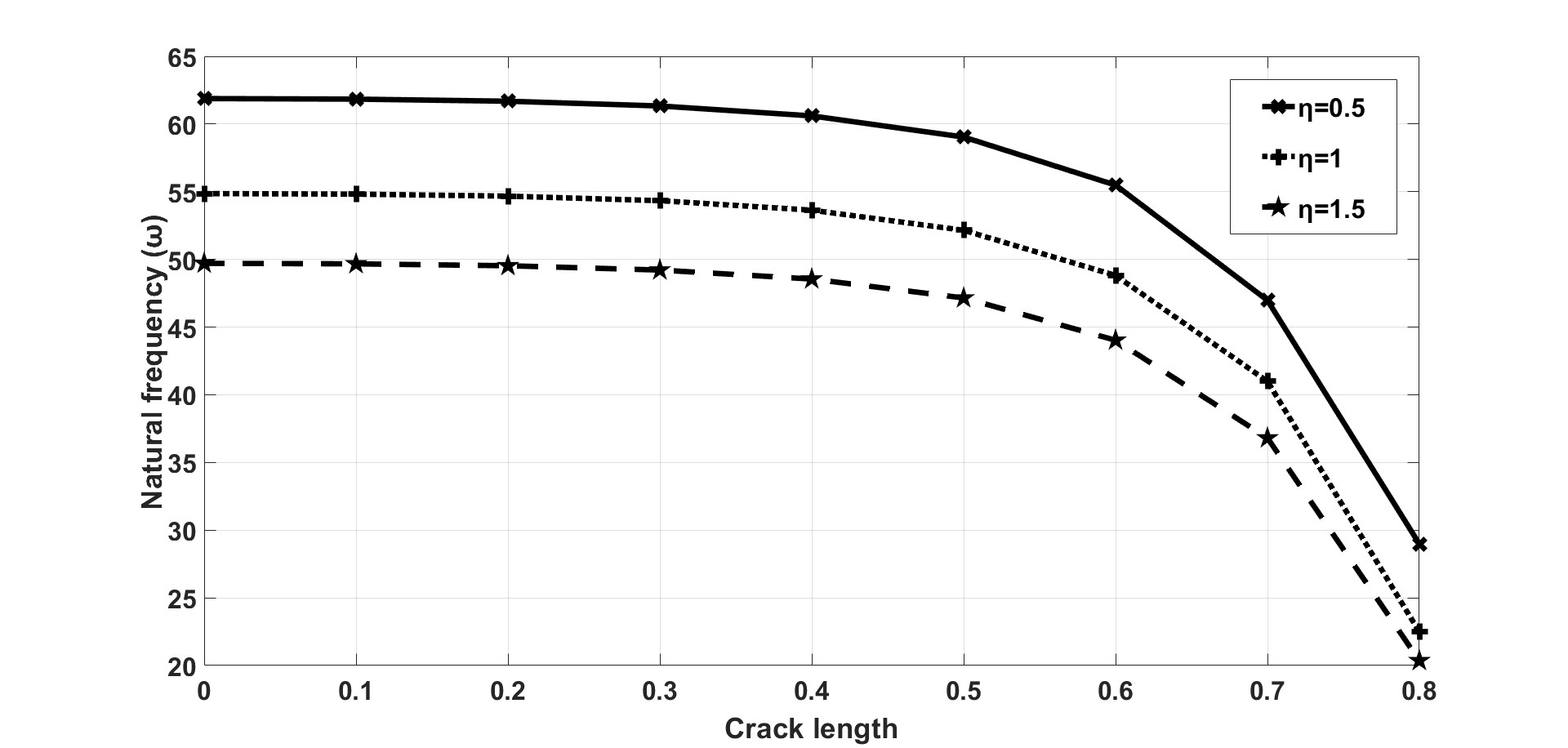}}
\caption{Natural frequency versus crack length for $\alpha=0.8$.}
\label{Fig.6}
\centering
\end{figure}

\begin{figure}[H]
\centerline{\includegraphics[width=14cm,height=8cm]{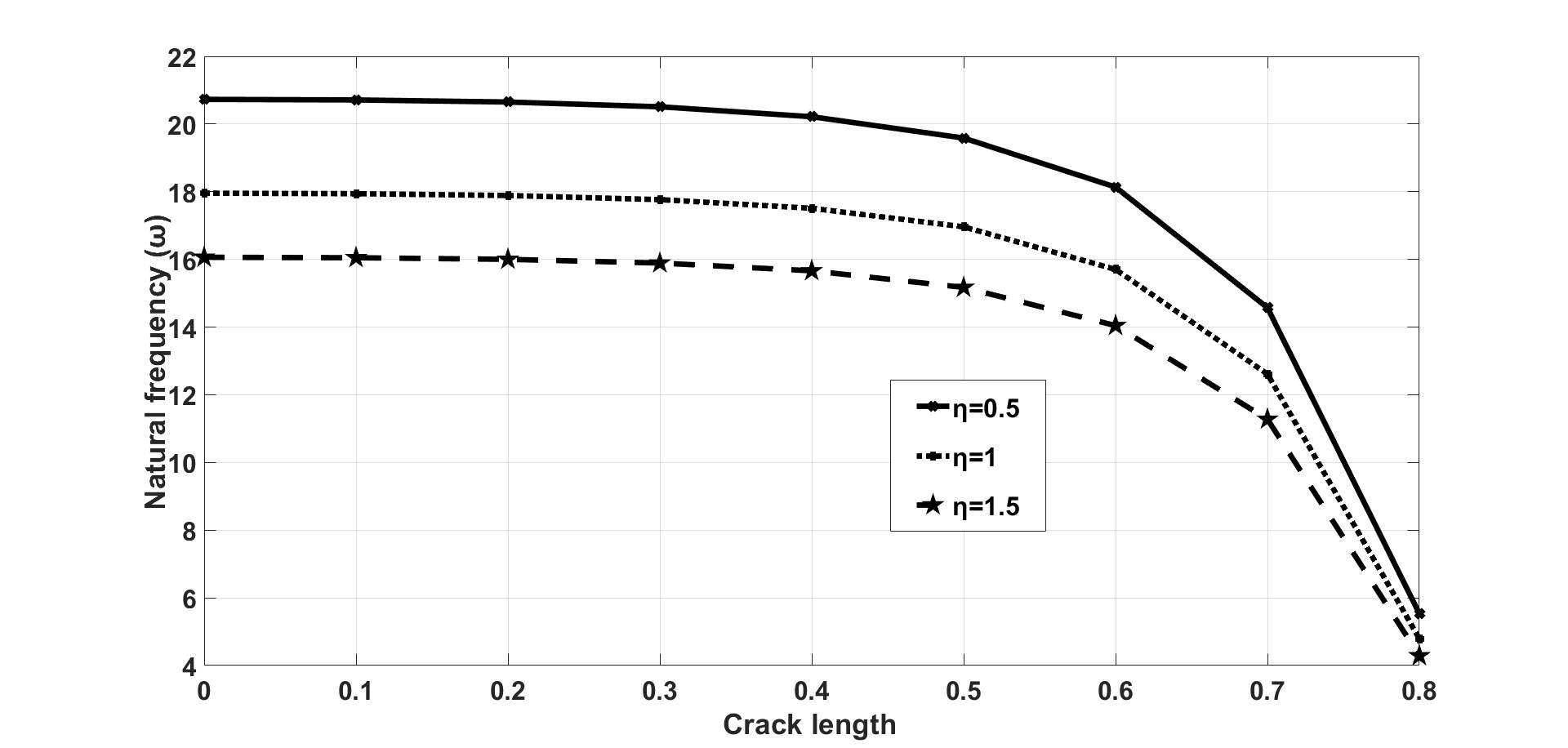}}
\caption{Natural frequency versus crack length for $\alpha=0.6$.}
\label{Fig.7}
\centering
\end{figure}

\begin{figure}[H]
\centerline{\includegraphics[width=14cm,height=8cm]{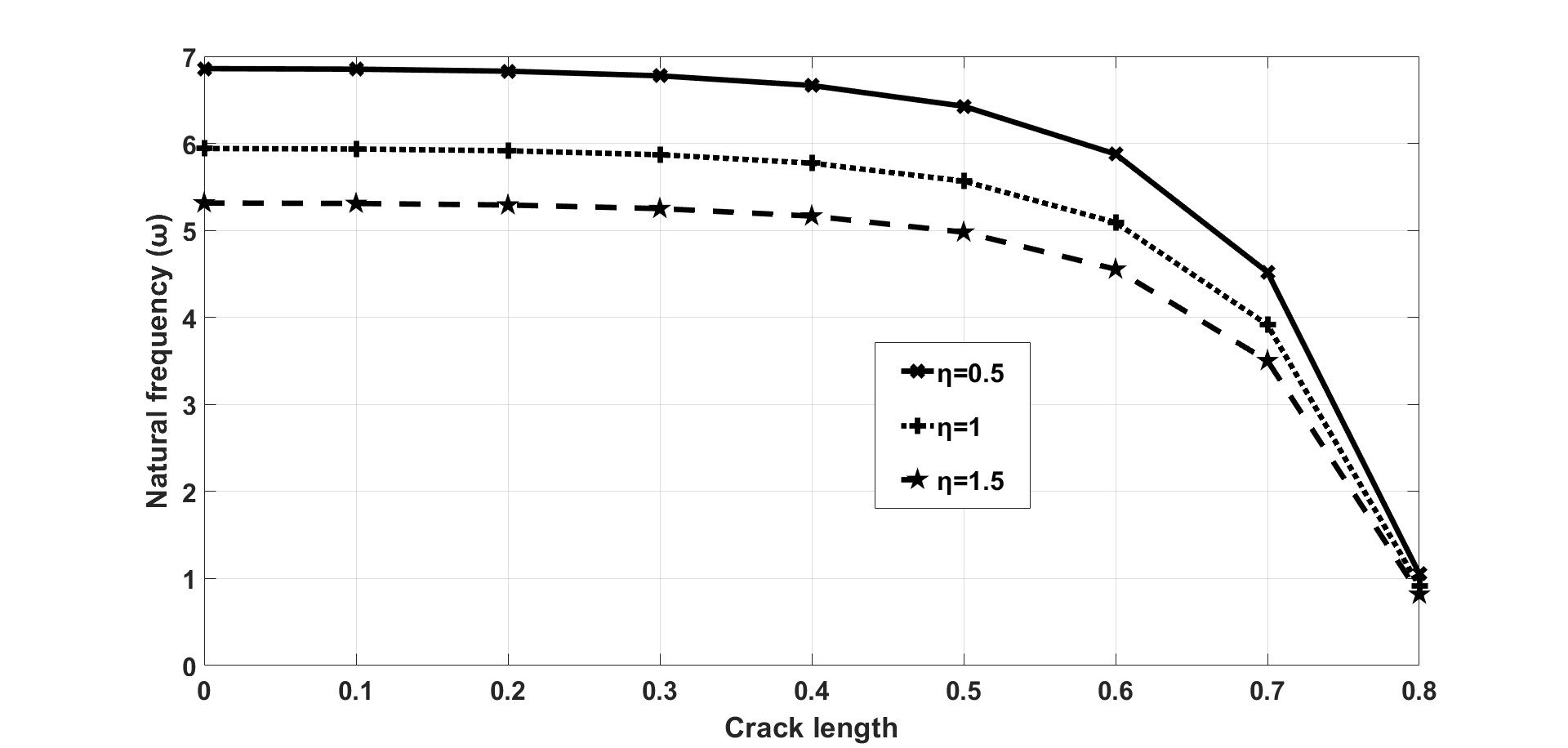}}
\caption{Natural frequency versus crack length for $\alpha=0.4$.}
\label{Fig.8}
\centering
\end{figure}

\begin{figure}[H]
\centerline{\includegraphics[width=14cm,height=8cm]{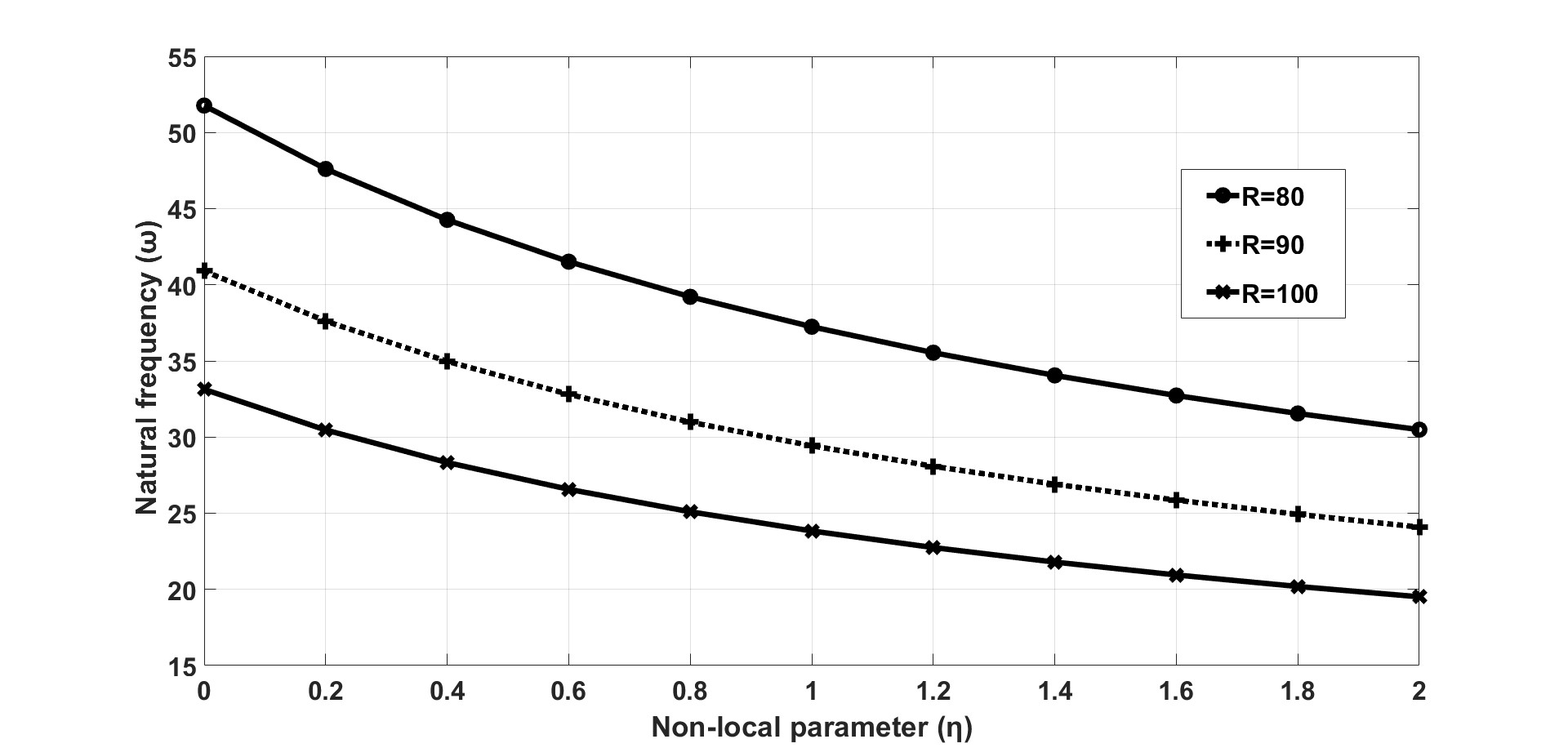}}
\caption{Natural frequency versus nonlocal parameter.}
\label{Fig.9}
\centering
\end{figure}

\begin{figure}[H]
\centerline{\includegraphics[width=14cm,height=8cm]{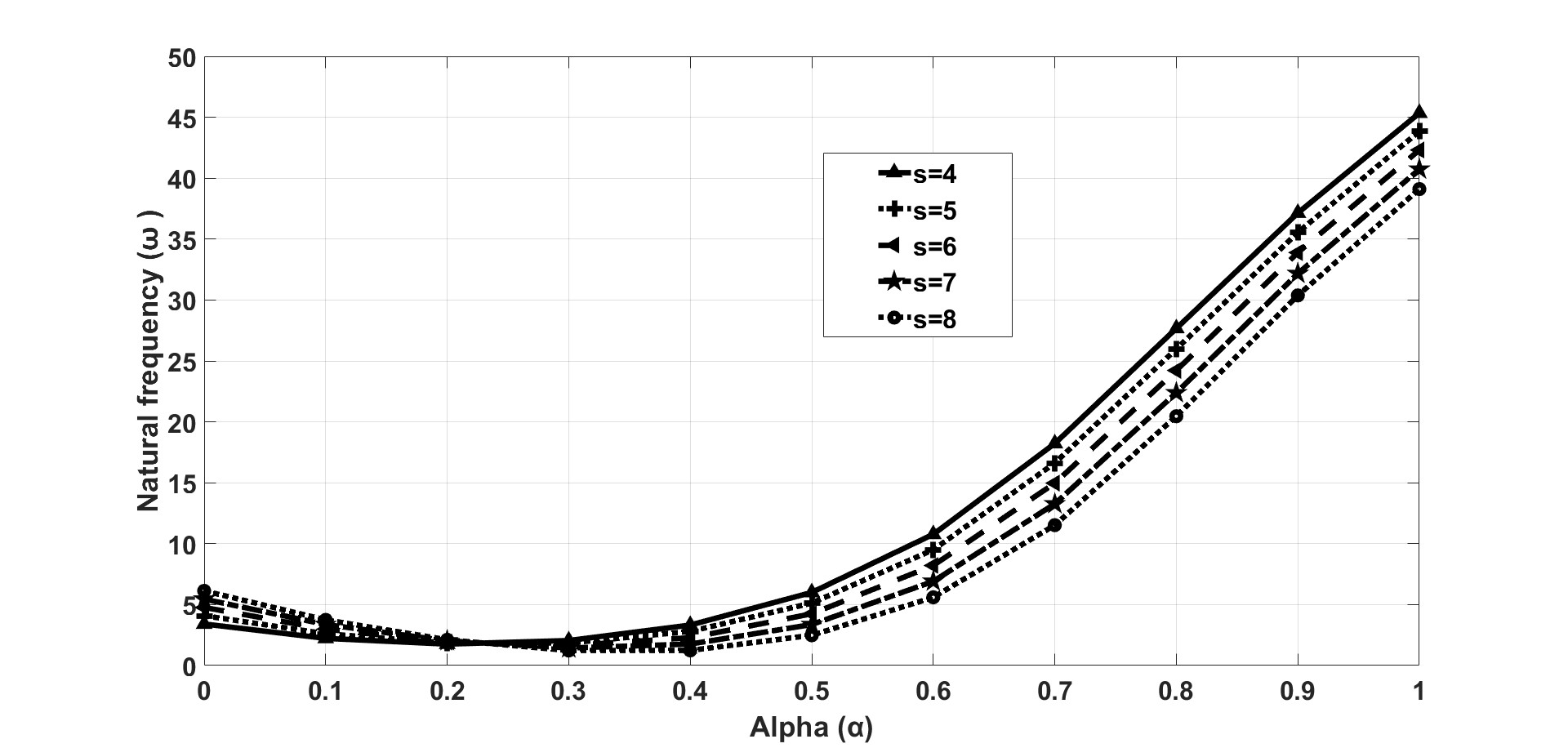}}
\caption{Natural frequency versus defect location.}
\label{Fig.10}
\centering
\end{figure}

\begin{figure}[H]
\centerline{\includegraphics[width=14cm,height=8cm]{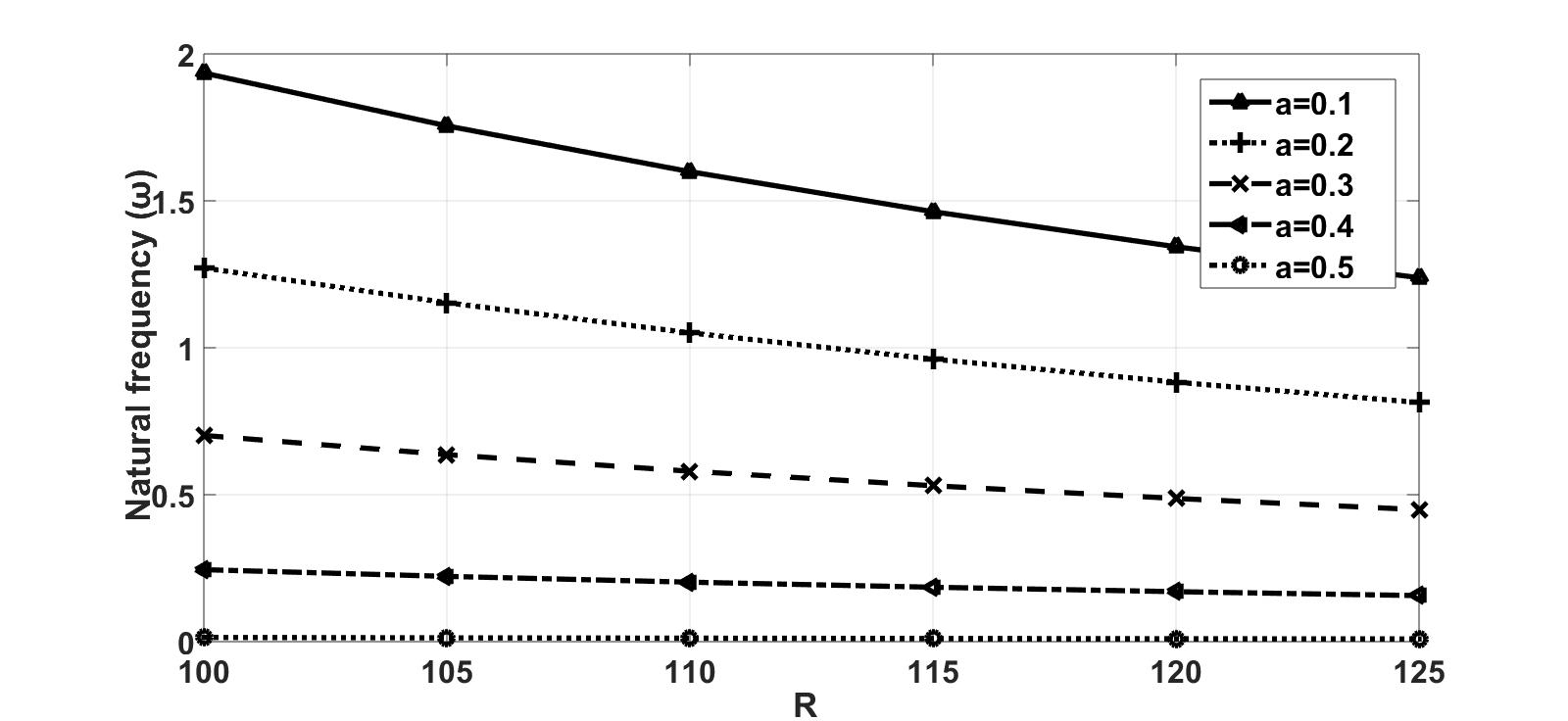}}
\caption{Natural frequency versus radius of the nano-arch.}
\label{Fig.11}
\centering
\end{figure}

\begin{figure}[H]
\centerline{\includegraphics[width=14cm,height=8cm]{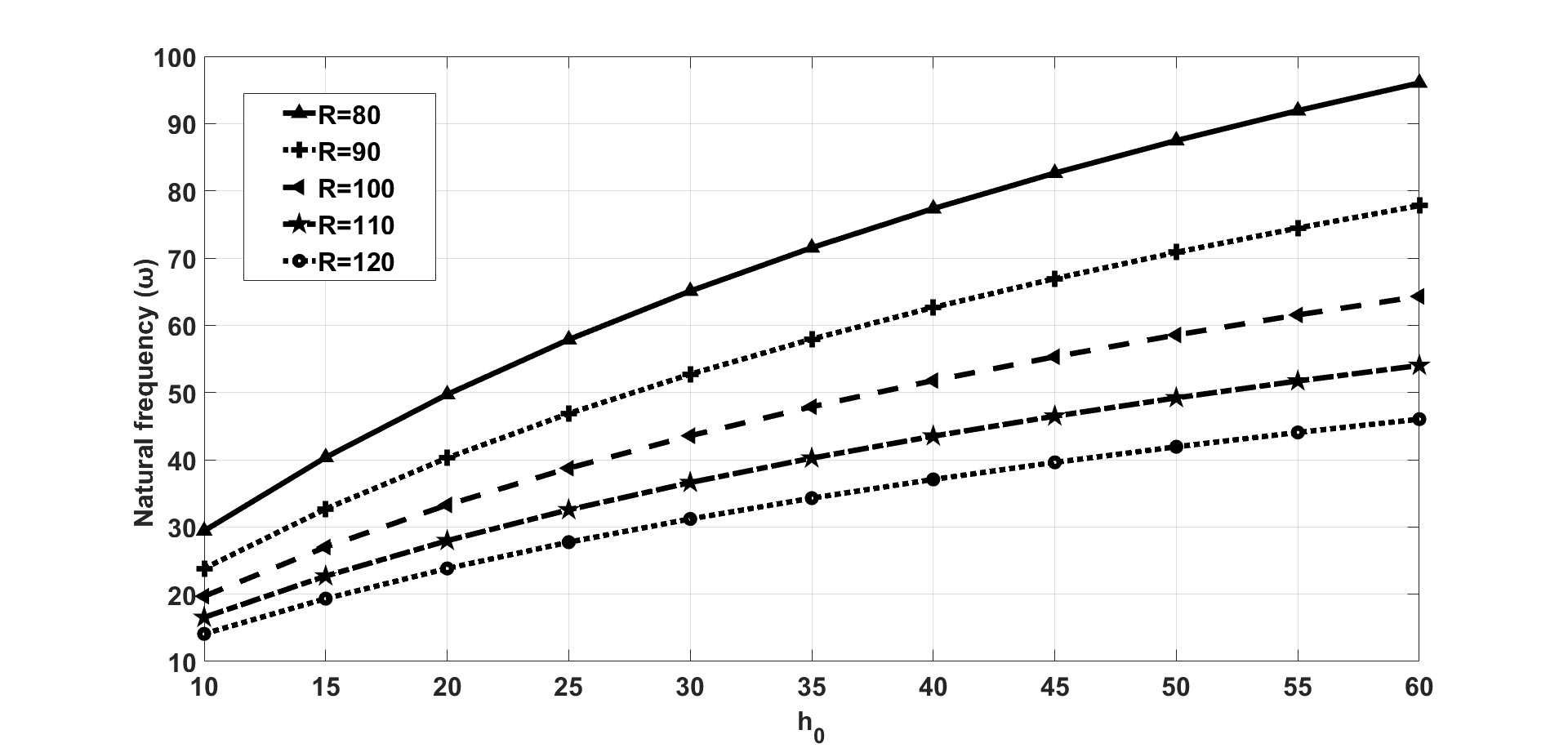}}
\caption{Natural frequency versus thickness $h_{0}$ of the nano-arch.}
\label{Fig.12}
\centering
\end{figure}

\begin{figure}[H]
\centerline{\includegraphics[width=14cm,height=8cm]{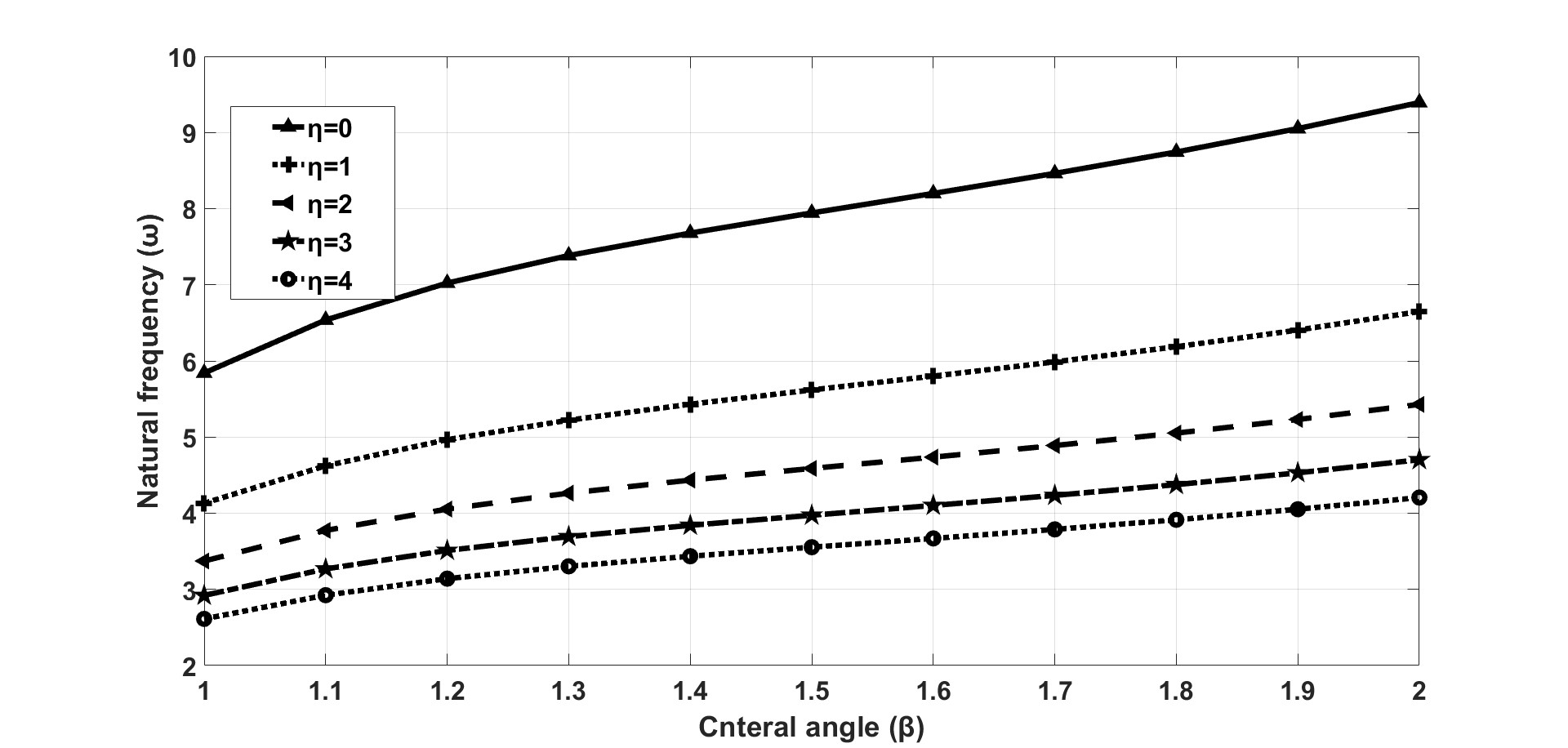}}
\caption{Natural frequency versus central angle of the nano-arch.}
\label{Fig.13}
\centering
\end{figure}

\begin{figure}[H]
\centerline{\includegraphics[width=14cm,height=8cm]{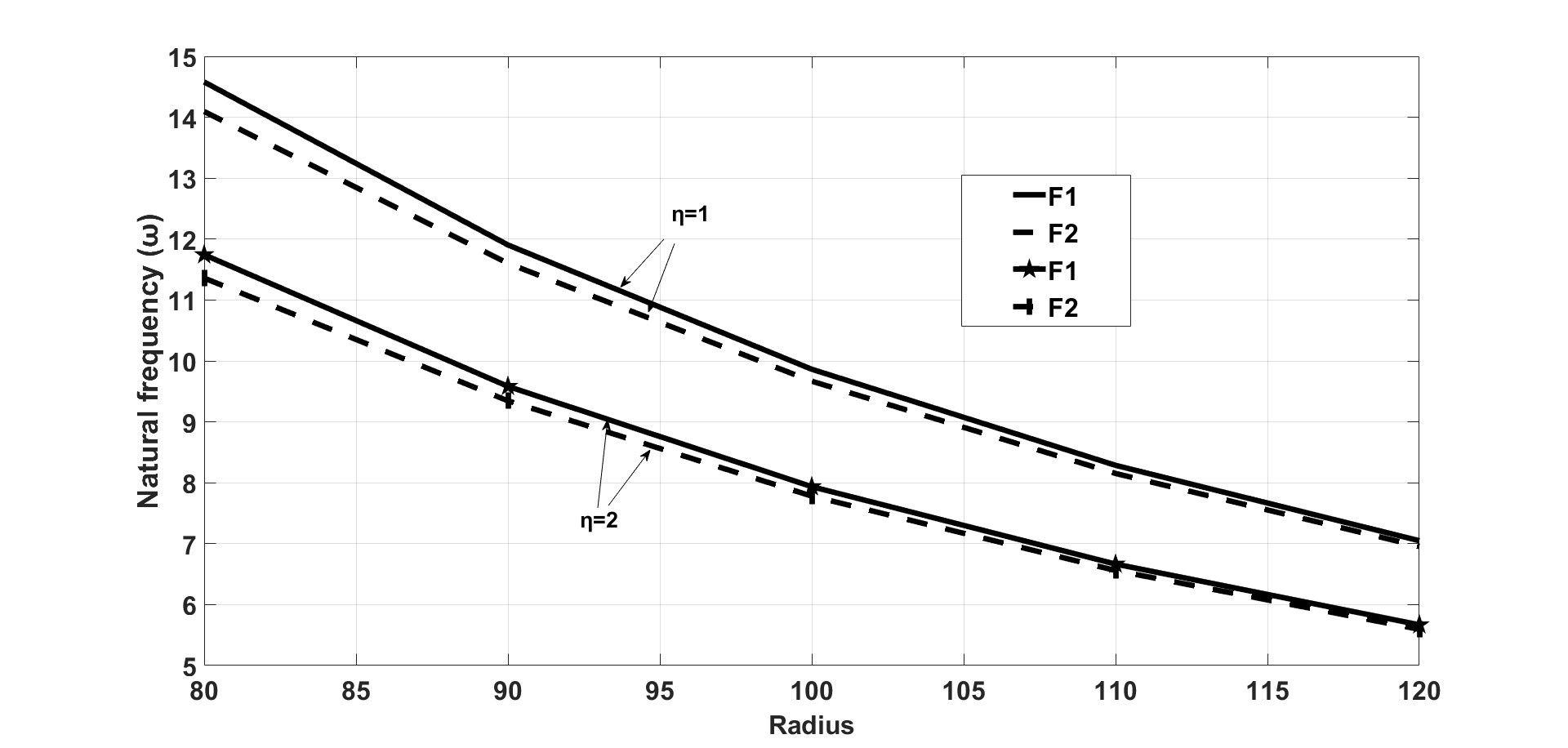}}
\caption{Natural frequency versus radius of the nano-arch for different shape functions.}
\label{Fig.14}
\centering
\end{figure}

Numerical results are obtained for a specimen made of a nanomaterial with the modulus of elasticity $E=7\times10^{11}\;Pa$, $\nu=0.3$, $\rho=10\;kg/m^3$. The radius of the middle line of the arch $R = 110\;nm$ and the width $b=1\;nm$ if the text does not contain any numerical evaluations of these quantities. The results of the calculations are presented in Figs. \ref{Fig.3}-\ref{Fig.14} and Tables \ref{tab:table1} and \ref{tab:table2}. In Figs. \ref{Fig.3}-\ref{Fig.5}, the natural frequency of the nano-arch is depicted versus the thickness of the nano-arch
for different values of the radius. Here $R=80-120\;nm$ and Figs. \ref{Fig.3}-\ref{Fig.5} correspond to $\alpha=0.8, 0.6$ and $0.4$, respectively. Looking at the figures, one can draw a conclusion that the smaller the radius, the higher the natural frequency of the nano-arch is for each given value of the thickness. Also, one can see from the figures that thicker nano-arches correspond to higher natural frequency values, provided the radius $R$ is fixed. The figures reveal another noteworthy conclusion: as the value of $\alpha$ increases, the impact of the thickness on the natural frequency of the nano-arches becomes greater.
Three different curves depicting the relationship between natural frequencies and crack length are presented in Figs. \ref{Fig.6}-\ref{Fig.8}, each curve corresponds to different values of the nonlocal parameter $\eta$. Specifically, Figs. \ref{Fig.6}-\ref{Fig.8} correspond to the values of $\alpha= 0.8, 0.6, 0.4$, respectively. The figures indicate that the natural frequency decreases as the crack length increases. Notably, a significant decrease in frequency occurs when the crack length surpasses 0.5. Additionally, the figures reveal that lower values of alpha are associated with decreased natural frequencies.
Fig. \ref{Fig.9} illustrates the relationship between the natural frequency of the nano-arch and the nonlocal parameter. The curves on the graph correspond to different values of the radius. It is evident from the figure that as the nonlocal parameter $\eta$ increases, the frequency of the nano-arch decreases. This behaviour is also observed when examining the natural frequency concerning the radius.
The relationship between the natural frequency and the location of the defect is depicted in Fig. \ref{Fig.10}. Each curve on the graph corresponds to a different crack length. It is noticeable from the figure that initially, the natural frequency decreases, but then it begins to increase. Once the crack location exceeds 0.5, the natural frequency experiences a significant increase. The effects of varying crack lengths can also be seen in the figure. 
In Fig. \ref{Fig.11}, the natural frequency is plotted against the radius of the nanoarch. Different curves in the figure correspond to the various positions of the defect. It can be seen in the figure that the natural frequency decreases with an increase in the radius. However, the influence of the radius on natural frequency is less observable for higher values of alpha.
Fig. \ref{Fig.12} illustrates the relationship between the natural frequency of the nanoarches and their thickness, with the thickness of the nanoarch varying between 0 and $\alpha$. The figure clearly shows that as the thickness of the nanoarch increases, so does the natural frequency. Different curves in Fig. \ref{Fig.12} represent various values of the radius, revealing the matter that the frequency decreases as the radius of the arch increases.

Fig. \ref{Fig.13} shows the relationship between the natural frequency and the central angle $\beta$ of the nano-arch. Different curves in Fig. \ref{Fig.13} correspond to different values of the nonlocal parameter $\eta$. It can be seen from Fig. \ref{Fig.13} that the natural frequency is higher for larger values of the central angle $\beta$. On the other hand, the smaller is the nonlocal parameter $\eta$ the higher is the natural frequency. In Fig. \ref{Fig.14}, the natural frequency is depicted against the radius of the nano-arch, showcasing various shape functions. From the figure, it is evident that the choice of shape functions significantly impacts the natural frequency of the arches.
Natural frequencies of nano-arches clamped at the left edge and free at the right-hand edge  are presented in Table \ref{tab:table1}. Table \ref{tab:table1}, corresponds to the arches with $\beta=30^0$. Table \ref{tab:table2} accommodates corresponding values of the natural frequency for arches fully clamped at both edges for different values of the parameter $\eta$ (here $0\leq \eta \leq 4$).
In Table \ref{tab:table1}, the results of the current work are compared with those obtained by Ganapathi et al. \cite{ganapathi2018vibration}, making use of the finite element method. The comparison of the eigenfrequencies found by different methods shows that the current method leads to somewhat overestimated values of natural frequencies. Nevertheless, the eigenfrequencies calculated by the current method are higher than those corresponding to Ganapathi et al. \cite{ganapathi2018vibration}.
The present calculation method is also applied to full nano-rings ($\beta=2\pi$). The obtained results of calculations are accommodated together with the results by Wang and Duan \cite{wang2008free} in Table \ref{tab:table2} for different values of the nonlocal parameter $\eta$. It can be seen from Table \ref{tab:table2} that the results obtained by the current method are close to those obtained by Wang and Duan. The current approach leads to the eigenfrequencies which are lower than the reseults of Wang and Duan \cite{wang2008free}.

\begin{table}[ht]
\centering
\caption{Natural frequency of the nano-arches for varying nonlocal parameter.}
\label{tab:table1}
\begin{tabular}{@{}cccc@{}} 
\toprule
 &        \multicolumn{2}{c}{$\beta=30^0$} & \\
\cline{1-4}
Mode & $\eta$ & Present & Ganapathi et al.\cite{ganapathi2018vibration} \\
\hline
1 & 0 & 5.1132 & 3.5078 \\
  & 1 & 4.8045 & 3.4289 \\
  & 2 & 4.5771 & 3.3546 \\
  & 3 & 4.1808 & 3.2847 \\
  & 4 & 3.8474 & 3.2187 \\
\cline{1-4}
2 & 0 & 21.6534 & 19.8974 \\
  & 1 & 18.9765 & 17.6058 \\
  & 2 & 17.0128 & 15.9352 \\
  & 3 & 15.2049 & 14.6692 \\
  & 4 & 14.7549 & 13.6734 \\
\cline{1-4}
\end{tabular}
\end{table}

\begin{table}[ht]
\caption{Natural frequency of circular nano-rings.}
\label{tab:table2}
{\begin{tabular}{@{}cccccc@{}}
\toprule
 & \multicolumn{2}{c}{$a=0.2$} &   &\multicolumn{2}{c}{$a=0.4$}    \\ 
\cline{2-3} \cline{5-6} 
Mode &  Present   & Wang \cite{wang2008free}                                                         
	&      &  Present  &          Wang \cite{wang2008free}                      \\    \hline
 1   & 5.6234   & 6.2       &    &
   3.2089   & 4.4         \\       
 2  &  41.2437   & 42.4        &    & 
  21.8362   & 23.6         \\         
 3&  127.2065   & 129.1         &     & 
   57.5028    & 59.5          \\
  4 &  274.3769   & 276.9        &   & 
  107.7692    & 110.8          \\
 5  &  485.1752   & 488.5        &   & 
  173.1745    & 176.3         \\    
    &     & Nano-ring with a defect   &      &  
      &         \\
  1  & 4.9908    & 5.9       &   & 
    3.2543  & 4.2        \\
 2  & 39.1549    & 40.8        &   &   
    21.8654  & 23.1         \\
  3  & 123.7280    & 125.0       &   &    
   56.5437   & 58.6         \\
 4  & 267.8756    & 270.0        &   &     
    107.0654   & 109.6        \\
5  & 475.5538    & 478.6         &  &     
    172.1767   & 174.9        \\
  \cline{1-6}                          
\end{tabular}}
\end{table}

\newpage

\bibliographystyle{amsplain}
\bibliography{main}
\end{document}